\documentclass[12pt]{amsart}
\usepackage{amssymb}
\usepackage{amsfonts}
\usepackage{color,soul}
\usepackage{appendix}

\newcommand{\nc}{\newcommand}
\nc{\thref}[1]{Theorem~\ref{theo:#1}}
\nc{\selabel}[1]{\label{sect:#1}}
\nc{\seref}[1]{Section~\ref{sect:#1}}
\nc{\lelabel}[1]{\label{lemm:#1}}
\nc{\leref}[1]{Lemma~\ref{lemm:#1}}
\nc{\prlabel}[1]{\label{prop:#1}}
\nc{\prref}[1]{Proposition~\ref{prop:#1}}
\nc{\colabel}[1]{\label{coro:#1}}
\nc{\coref}[1]{Corollary~\ref{coro:#1}}
\nc{\exlabel}[1]{\label{exam:#1}}
\nc{\exref}[1]{Example~\ref{exam:#1}}
\nc{\delabel}[1]{\label{defi:#1}}
\nc{\deref}[1]{Definition~\ref{defi:#1}}
\nc{\eqlabel}[1]{\label{equa:#1}}
\nc{\relabel}[1]{\label{rema:#1}}
\nc{\reref}[1]{Lemma~\ref{rema:#1}}
\providecommand{\operatorname}[1]{\mathrm{#1}\,}
\nc{\Hom}{\operatorname{Hom}} \nc{\Mor}{\operatorname{Mor}}
\nc{\Aut}{\operatorname{Aut}} \nc{\Ann}{\operatorname{Ann}}
\nc{\Ker}{\operatorname{Ker}} \nc{\Trace}{\operatorname{Trace}}
\nc{\Char}{\operatorname{Char}} \nc{\Mod}{\operatorname{Mod}}
\nc{\End}{\operatorname{End}} \nc{\Spec}{\operatorname{Spec}}
\nc{\Span}{\operatorname{Span}} \nc{\sgn}{\operatorname{sgn}}
\nc{\Id}{\operatorname{Id}} \nc{\Com}{\operatorname{Com}}
\nc{\rank}{\operatorname{rank}}
\nc{\Clausen}{\operatorname{Cl}}
\nc{\Li}{\operatorname{Li}}
\nc{\Ls}{\operatorname{Ls}}

\let\:=\colon

\newtheorem{de}{Definition}[section]
\newtheorem{lm}[de]{Lemma}
\newtheorem{pr}[de]{Proposition}
\newtheorem{co}[de]{Corollary}
\newtheorem{re}[de]{Remark}
\newtheorem{res}[de]{Remarks}
\newtheorem{te}[de]{Theorem}
\newtheorem{ex}[de]{Example}
\newtheorem{exs}[de]{Examples}

\def\bex{\begin{ex}}
\def\eex{\end{ex}}
\def\bexs{\begin{exs}}
\def\eexs{\end{exs}}
\def\bl{\begin{lm}}
\def\el{\end{lm}}
\def\bc{\begin{co}}
\def\ec{\end{co}}
\def\bt{\begin{te}}
\def\et{\end{te}}
\def\bpr{\begin{pr}}
\def\epr{\end{pr}}
\def\br{\begin{re}}
\def\er{\end{re}}
\def\brs{\begin{res}}
\def\ers{\end{res}}
\def\bd{\begin{de}}
\def\ed{\end{de}}
\def\be{\begin{equation}}
\def\ee{\end{equation}}
\def\bea{\begin{eqnarray*}}
\def\eea{\end{eqnarray*}}
\def\bp{\begin{proof}}
\def\ep{\end{proof}}

\def\qed{\hfill\Box}

\let\:=\colon

\begin{document}

\title[\hfil Zagier formula for multiple zeta values]
{Elementary proofs of Zagier's formula for multiple zeta values and its odd variant}

\begin{abstract}
In this paper, we give elementary proofs of Zagier's formula for multiple zeta values involving Hoffman elements and its odd variant due to Murakami. Zagier's formula was a key ingredient in the proof of Hoffman's conjecture. Moreover, using the same approach, we prove Murakami's formula for multiple $t$-values. This formula is essential in proving a Brown type result which asserts that each multiple zeta value is a $\mathbb{Q}$-linear combination of multiple $t$-values of the same weight involving $2$'s and $3$'s.       
\end{abstract}

\author{Li Lai, Cezar Lupu, Derek Orr}



\maketitle

\section{Introduction and statement of results}

Multiple zeta values (abbreviated MZV's) of weight $|k|=k_{1}+\ldots +k_{r}$ and depth (length) $r$, defined by 

$$\displaystyle\zeta(k_{1}, k_{2}, \ldots, k_{r})=\sum_{1\leq n_{1}<n_{2}<\ldots <n_{r}}\frac{1}{n_{1}^{k_{1}}n_{2}^{k_{2}}\ldots n_{r}^{k_{r}}}$$
for positive integers $k_{1}, k_{2}, \ldots, k_{r}$ with $k_{r}>1$, have been a subject of great interests over the past 30 years. These numbers are intimately connected to various objects in mathematics and physics such as knot invariants, quantum groups, Galois representations of the fundamental group of $\mathbb{P}^{1}-\{0, 1, \infty\}$, mixed Tate motives or Feynman amplitudes in perturbative quantum field theory. 
\vspace{0.3cm}

 Starting with Euler's remarkable identity $\zeta(1, 2)=\zeta(3)$, MZV's generate an algebra and they satisfy a variety of relations. A challenging problem concerning MZV's is to determine all $\mathbb{Q}$-linear and polynomial relations. In this direction, M. Hoffman \cite{Hoffman1} conjectured that all multiple zeta values of a given weight can be expressed as a $\mathbb{Q}$-linear combination of MZV's of the same weight involving $2$'s and $3$'s. 
 For example, in weight $5$ all MZV's are linear combinations of $\zeta(2, 3)$ and $\zeta(3, 2)$, and in weight $7$, all MZV's are linear combinations of $\zeta(2, 2, 3)$, $\zeta(2, 3, 2)$ and $\zeta(3, 2, 2)$. 
\vspace{0.3cm}

Hoffman's conjecture was proved by F. Brown \cite{Brown} in 2012 using motivic multiple zeta values. One of the key ingredients in his proof was Zagier's formula \cite{Zagier1} for a specific family of MZV's which we call nowadays the Hoffman family, 

$$\displaystyle H(a, b)=\zeta(\underbrace{2, 2, \ldots, 2}_{\text{$a$}}, 3, \underbrace{2, 2, \ldots, 2}_{\text{$b$}}).$$

Moreover, Brown showed that $H(a, b)$ can be written as a $\mathbb{Q}$-linear combination of products $\pi^{2m}\zeta(2n+1)$ with $m+n=a+b+1$. However, Brown's motivic formalism does not provide the precise value of those coefficients from the linear combination of products mentioned above. This task was accomplished by Zagier \cite{Zagier1} where he provided the formula for $H(a, b)$ and its proof. This formula of Zagier sparked a lot of interests in the mathematical community due to its wide range and unexpected applications. For example, one such application was proving Broadhurst's zig-zag conjecture by Brown and Schnetz \cite{Brown-Schnetz}. More details about the formula are given in the next section.
\vspace{0.3cm}

First, let us recall the much easier formula for the simplest Hoffman family of MZV's \cite{Hoffman, Zagier, Lupu}, namely

$$\zeta(\underbrace{2, 2, \ldots, 2}_{\text{$n$}})=\frac{\pi^{2n}}{(2n+1)!}.$$

This Hoffman-Zagier evaluation from above seems to be present in many other evaluations which concern different families of multiple zeta values. One such example is Zagier's formula \cite{Zagier1} for Hoffman elements, which reads as
\vspace{0.3cm}

\bt \label{Theorem_Zagier}
Let $\displaystyle H(n):=\zeta(\underbrace{2, 2, \ldots, 2}_{\text{$n$}})$. 
We have 
\begin{equation*}\displaystyle H(a, b):=\zeta(\underbrace{2, 2, \ldots, 2}_{\text{$a$}}, 3, \underbrace{2, 2, \ldots, 2}_{\text{$b$}} )=2\sum_{k=1}^{a+b+1}(-1)^kc_{a, b}^{k}\zeta(2k+1)H(a+b+1-k),
\end{equation*}
where $\displaystyle c_{a, b}^{k}=\binom{2k}{2a+2}-\left(1-\frac{1}{2^{2k}}\right)\binom{2k}{2b+1}$ and $a, b\geq 0$ are integers. 
\et
\vspace{0.4cm}
 
In a similar way, M. Hoffman \cite{Hoffman2} defined the multiple $t$-values (``odd variant'' of MZV's),

$$\displaystyle t(k_{1}, k_{2}, \ldots, k_{r})=2^{-|k|}\zeta(k_{1}, k_{2}, \ldots, k_{r}; -\frac{1}{2}, -\frac{1}{2}, \ldots, -\frac{1}{2})=$$
$$\displaystyle=\sum_{1\leq n_{1}<n_{2}<\ldots<n_{r}}\frac{1}{(2n_{1}-1)^{k_{1}}(2n_{2}-1)^{k_{2}}\ldots (2n_{r}-1)^{k_{r}}}$$
$$\displaystyle=\sum_{1\leq n_{1}<n_{2}<\ldots <n_{r}:\\ n_{i} \operatorname{odd}}\frac{1}{n_{1}^{k_{1}}n_{2}^{k_{2}}\ldots n_{r}^{k_{r}}},$$
of weight $|k|=k_{1}+\ldots+k_{r}$ and depth $r$. Here $\displaystyle\zeta(k_{1}, k_{2}, \ldots, k_{r}; a_{1}, a_{2}, \ldots, a_{r})$ stands for the multiple Hurwitz zeta values,
$$\displaystyle\zeta(k_{1}, k_{2}, \ldots, k_{r}; a_{1}, a_{2}, \ldots, a_{r})=\displaystyle\sum_{1\leq n_{1}<n_{2}<\ldots <n_{r}}\frac{1}{(n_{1}+a_{1})^{k_{1}}(n_{2}+a_{2})^{k_{2}}\ldots (n_{r}+a_{r})^{k_{r}}}.$$
\vspace{0.3cm}

As it has been already highlighted in \cite{Hoffman2} multiple $t$-values have remarkable parallels to and contrasts with MZV's. For example, the analogue of Euler's identity is given by $\displaystyle t(1, 2)=-\frac{1}{2}t(3)+t(2)\log2$, where $t(2)=\frac{\pi^2}{8}$. This tells us that even though multiple $t$-values generate an algebra, this is quite different in some ways from the algebra of MZV's. More details in this direction are presented in \cite{Hoffman2} and other related results can be found in \cite{Chung, Murty-Sinha, Shen-Jia, Zhao1}. 
\vspace{0.3cm}

In a slightly different direction, Murakami \cite{Murakami} proved that in some cases multiple $t$-values are linear combinations of multiple zeta values. More exactly, in the spirit of Brown's theorem \cite{Brown} (Hoffman's conjecture), Murakami shows that every multiple zeta value is a $\mathbb{Q}$-linear combination of elements $\{t(k_{1}, \ldots, k_{r}): k_{1}, \ldots, k_{r}\in \{2, 3\}\}$. Again, an evaluation of Zagier-type is needed for the Hoffman family of multiple $t$-values,  
$$\displaystyle T(a, b)=t(\underbrace{2, 2, \ldots, 2}_{\text{$a$}}, 3, \underbrace{2, 2, \ldots, 2}_{\text{$b$}}).$$

As we stated above, first let us recall the the analogue formula for the simplest Hoffman elements \cite{Chung, Lupu, Murty-Sinha, Shen-Jia, Zhao1},

$$t(\underbrace{2, 2, \ldots, 2}_{\text{$n$}})=\frac{\pi^{2n}}{2^{2n}(2n)!}.$$

Also, the analogue for multiple $t$-values of Zagier's formula (Theorem \ref{Theorem_Zagier}) is given by
\vspace{0.3cm}

\bt \label{Theorem_odd_variant}
Let $\displaystyle T(n):=t(\underbrace{2, 2, \ldots, 2}_{\text{$n$}})$. Then we have

\begin{equation*}
    T(a,b) = \sum_{k=1}^{a+b+1} (-1)^{k+1}\left(\binom{2k}{2a+1}+\binom{2k}{2b+1}\left(1-\frac{1}{2^{2k}}\right)\right)\frac{1}{2^{2k}}T(a+b+1-k)\zeta(2k+1).
\end{equation*}

\et
\vspace{0.4cm}

In this paper, building upon previous work from \cite{Lupu} we derive elementary and direct proofs for both $H(a, b)$ (Theorem \ref{Theorem_Zagier}) and $T(a, b)$ (Theorem \ref{Theorem_odd_variant}) using the same approach with identities explored in \cite{Orr}. Our proofs are significantly simpler than those given in \cite{Murakami, Zagier1} and rely on the Taylor series of integer powers of arcsin. The central point of this paper is the parallel between the two formulas for $H(a, b)$ and $T(a, b)$. More about this is explained in the next sections. Furthermore, we explore the method used in the proof of the two formulas to obtain some arithmetic information about $\zeta(2k+1)/\pi^{2k+1}$.
\vspace{0.3cm}

\textbf{Organization of the paper and description of the proofs.} The paper is organized as follows. In section \ref{Preliminaries}, we collect some basic tools about Clausen functions (integrals), cotangent integrals and their relations with Riemann zeta values. In section \ref{Section_Zagier_formula} we present the proof of Zagier's formula for $H(a, b)$. The proof will consist of two major theorems. First, using the Taylor series expansion of $\arcsin^{2r}(x)$, we can express $H(a, b)$ in terms of an integral of the type $\int_0^{\frac{\pi}{2}}x^{2a+2}(1-\frac{2}{\pi}x)^{2b+1}\cot xdx$. Second, we express the previous mentioned integral as a $\mathbb{Q}$-linear combination of powers of $\pi$ and odd zeta values which will give us the exact formula. Theorem \ref{Theorem_Zagier} will follow from a combination of Theorem \ref{Theorem_H(a,b)} and Theorem \ref{Theorem_hat_H(a,b)}. In section \ref{Section_odd_variant}, we derive the proof of Zagier-type formula for $T(a, b)$ using the similar tools as in section \ref{Section_Zagier_formula}. This time, the difference is that we will rely on the Taylor series of $\arcsin^{2r+1}(x)$ and $\int_0^{\frac{\pi}{2}}x^{2a+1}(1-\frac{2}{\pi}x)^{2b+1}\cot xdx$. Theorem \ref{Theorem_odd_variant} will follow from Theorems \ref{Theorem_T(a,b)} and \ref{Theorem_hat_T(a,b)}. In the last section, we will give a brief account on how Zagier's formula leads to some non-trivial arithmetic properties of $\frac{\zeta(2k+1)}{\pi^{2k+1}}$.    
\vspace{0.3cm}

\textbf{Acknowledgements.} We would like to thank Vlad Matei and Cosmin Pohoa\c t\u a for their careful reading of the manuscript which led to a better presentation of this paper. Also, we thank Florian Luca, Wadim Zudilin for illuminating discussions over an earlier version of this manuscript and to Masanobu Kaneko for his interests in this paper.  

\section{Preliminaries}\label{Preliminaries}

The Clausen function (integral) and the higher order Clausen functions are defined by
$$\displaystyle\operatorname{Cl}_{2}(\theta)=-\int_0^{\theta}\log\left(2\sin\frac{t}{2}\right)dt=\sum_{k=1}^{\infty}\frac{\sin(k\theta)}{k^2},$$
$$\displaystyle\operatorname{Cl}_{2m}(\theta)=\sum_{k=1}^{\infty}\frac{\sin(k\theta)}{k^{2m+1}}, \operatorname{Cl}_{2m+1}(\theta)=\sum_{k=1}^{\infty}\frac{\cos(k\theta)}{k^{2m+1}}.$$

Some important particular values include
$$\displaystyle\operatorname{Cl}_{2m}(\pi)=0, \operatorname{Cl}_{2m+1}(\pi)=-\frac{(4^m-1)\zeta(2m+1)}{4^m} ,$$
and $$\displaystyle\operatorname{Cl}_{2m}\left(\frac{\pi}{2}\right)=\beta(2m), \operatorname{Cl}_{2m+1}\left(\frac{\pi}{2}\right)=-\frac{(4^m-1)\zeta(2m+1)}{2^{4m+1}},$$
where $\beta(s)=\sum_{n=0}^{\infty}\frac{(-1)^n}{(2n+1)^s}, \operatorname{Re} s>0$ is the Dirichlet beta function.
\vspace{0.3cm}

As we have already highlighted in \cite{Lupu1, Orr}, we can express the cotangent integral in terms of Clausen functions.

\bl[\cite{Lupu1, Orr}]\label{Lemma_cot_integral}
For $p \in \mathbb{N}$ and $|z|<1$, we have

\begin{align*}
\int_{0}^{\pi z} x^{p} \cot (x)  dx=&(\pi z)^{p} \sum_{k=0}^{p} \binom{p}{k} \frac{k!(-1)^{\lfloor(k+3) / 2\rfloor}}{(2 \pi z)^{k}} \operatorname{Cl}_{k+1}(2 \pi z) \\
&+\delta_{\lfloor p / 2\rfloor, p / 2} \frac{p !(-1)^{p / 2}}{2^{p}} \zeta(p+1).
\end{align*}
\el

Now, for $z=\frac{1}{2}$ the above lemma reads as

\bl
\begin{align*}
\int_{0}^{\pi / 2} x^{p} \cot (x) \mathrm{d} x=&\left(\frac{\pi}{2}\right)^{p}\left(\log 2+\sum_{k=1}^{\lfloor p / 2\rfloor} \frac{p !(-1)^{k}\left(4^{k}-1\right)}{(p-2 k) !(2 \pi)^{2 k}} \zeta(2 k+1)\right) \\
&+\delta_{\lfloor p / 2\rfloor, p / 2} \frac{p !(-1)^{p / 2} \zeta(p+1)}{2^{p}}.
\end{align*}
\el  
\vspace{0.3cm}

An alternative statement of above lemma is given by the following

\bl\label{Lemma_Orr}
For any polynomial $P(x) \in \mathbb{C}[x]$ with $P(0) = 0$, we have
\begin{align}
\int_{0}^{1} P(x) \cot\left( \frac{\pi x}{2} \right) {\rm d} x = 2P(1) \frac{\log 2}{\pi} + &2\sum_{k=1}^{ \lfloor (\deg P) /2 \rfloor} (-1)^{k} P^{(2k)}(1) \cdot \left( 1 -\frac{1}{2^{2k}} \right) \frac{\zeta(2k+1)}{\pi^{2k+1}} \notag \\
+&2\sum_{k=1}^{\lfloor (\deg P) /2 \rfloor} (-1)^{k} P^{(2k)}(0) \cdot \frac{\zeta(2k+1)}{\pi^{2k+1}}. \label{Lemma}
\end{align}
(Where $P^{(2k)}(x)$ stands for the $(2k)$-th derivative of $P(x)$.)
\el

\section{Zagier's formula for $H(a, b)$}\label{Section_Zagier_formula}

The proof given by Zagier \cite{Zagier1} of Theorem \ref{Theorem_Zagier} was indirect and it goes along the following lines. Consider the generating functions

$$\displaystyle F(x, y)=\sum_{a, b\geq 0}(-1)^{a+b+1}H(a, b)x^{2a+2}y^{2b+1},$$
and
$$\displaystyle\hat{F}(x, y)=\sum_{a, b\geq 0}(-1)^{a+b+1}\hat{H}(a, b)x^{2a+2}y^{2b+1},$$
where $\hat{H}(a, b)$ is the right-hand side quantity from Theorem \ref{Theorem_Zagier}. Furthermore, Zagier showed that $$\displaystyle\frac{\pi}{\sin(\pi y)}F(x, y)=\frac{\pi}{\sin(\pi y)}\hat{F}(x, y).$$
\vspace{0.3cm}

However, using transformation formulas for the hypergeometric function ${}_3F_{2}$, Li \cite{Li} simplified Zagier's proof. Other alternative proofs using different transformation formulas are given in \cite{Hessami-Pilehrood2, Lee-Peng}. 
\vspace{0.3cm}

Now, let us recall that the Taylor series for even integer powers of arcsin by comparing the coefficients of like powers of $\lambda$ in the formulas

$$\displaystyle\cos(\lambda\arcsin x)={}_{2} F_{1}\left(\frac{\lambda}{2}, -\frac{\lambda}{2}; \frac{1}{2}; x^2\right),$$
where ${}_2F_{1}$ is the Gauss hypergeometric function defined by $\displaystyle {}_{2}F_{1}(a,b;c;z) := \sum_{n=0}^{\infty} \frac{(a)_{n}(b)_{n}}{(c)_{n}}\frac{z^n}{n!}$
with $(p)_{0} = 1$ and $(p)_{n} = p(p+1)(p+2)\cdots(p+n-1)$ for $n>0$. This implies the following Taylor series expansion,

\begin{equation} \label{Equation_Taylor_arcsin_even_powers}
   \displaystyle\frac{\arcsin^{2r}(x)}{(2r)!}=\frac{1}{4^r}\sum_{n=1}^{\infty}\frac{4^n}{n^2\binom{2n}{n}}\cdot x^{2n}\cdot\sum_{n_{1}<n_{2}<\ldots n_{r-1}<n}\frac{1}{n_{1}^2n_{2}^2\ldots n_{r-1}^2}.
\end{equation}
\vspace{0.3cm}

Theorem \ref{Theorem_Zagier} will follow from a combination of Theorem \ref{Theorem_H(a,b)} and Theorem \ref{Theorem_hat_H(a,b)}. First, the following lemma tells us that the moments of $\arccos^{2b+1} x$ are related to the tail multiple harmonic numbers.

\bl \label{Lemma_arccos_odd_moments}
For any positive integer $n$ and nonnegative integer $b$, we have

\[ \int_{0}^{1} x^{2n-1}  \frac{(2\arccos x)^{2b+1}}{(2b+1)!} dx = \frac{\binom{2n}{n}}{n4^n}\frac{\pi}{2} \sum_{n<m_1<\cdots<m_b} \frac{1}{m_1^2 \cdots m_b^2}. \]
When $b=0$, the sum $\sum_{n<m_1<\cdots<m_b}$ is understood as $1$.
\el

\vspace{0.2cm}

\textit{Proof of Lemma \ref{Lemma_arccos_odd_moments}.} Changing the variable $x \mapsto \cos x$ and integrating by parts, we have
\begin{align*}
\int_{0}^{1} x^{2n-1}  \frac{(2\arccos x)^{2b+1}}{(2b+1)!} dx &=  \int_{0}^{\pi/2} \cos^{2n-1}(x) \sin(x)   \frac{(2x)^{2b+1}}{(2b+1)!} dx  \\
&= \frac{1}{n} \int_{0}^{\pi/2} \cos^{2n} (x)  \frac{(2x)^{2b}}{(2b)!} dx.
\end{align*}
We denote
\[ I_{n,b}:= \int_{0}^{\pi/2} \cos^{2n} (x)  \frac{(2x)^{2b}}{(2b)!} dx. \]
It remains to show that
\begin{equation}\label{Equation_Inb}
 I_{n,b} =  \frac{\binom{2n}{n}}{4^n}\frac{\pi}{2} \sum_{n<m_1<\cdots<m_b} \frac{1}{m_1^2 \cdots m_b^2}.
\end{equation}
We perform induction on $n+b$ to show that \eqref{Equation_Inb} holds for any nonnegative integers $n,b$. Note that for $b=0$ this is Wallis' formula, and for $n=0$ this is the well known fact $H(b):=\zeta(\underbrace{2, 2, \ldots, 2}_{\text{$b$}}) = \pi^{2b}/(2b+1)!$. In the sequel, we assume $n,b>0$ and the equality \eqref{Equation_Inb} holds for smaller $n+b$.
\vspace{0.2cm}
Using integration by parts again, we have
\begin{align*}
I_{n,b} &= \int_{0}^{\pi/2} \cos^{2n-1} (x)  \frac{(2x)^{2b}}{(2b)!} d (\sin x) \\
&= \int_{0}^{\pi/2} (2n-1) \cos^{2n-2} (x) \sin^2 (x)  \frac{(2x)^{2b}}{(2b)!} dx -   \int_{0}^{\pi/2} \cos^{2n-1}(x) \sin(x) \frac{2\cdot(2x)^{2b-1}}{(2b-1)!} dx \\
&= (2n-1)I_{n-1,b} - (2n-1) I_{n,b} -  \int_{0}^{\pi/2} \cos^{2n-1}(x) \sin(x) \frac{2\cdot(2x)^{2b-1}}{(2b-1)!} dx,
\end{align*}
so
\begin{equation}\label{Equation_Inb_1}
I_{n,b} = \frac{2n-1}{2n} I_{n-1,b} - \frac{1}{2n} \int_{0}^{\pi/2} \cos^{2n-1}(x) \sin(x) \frac{2\cdot(2x)^{2b-1}}{(2b-1)!} dx.
\end{equation}
Also, we have
\vspace{0.2cm}
\begin{align*}
&  \int_{0}^{\pi/2} \cos^{2n-1}(x) \sin(x)  \frac{2\cdot(2x)^{2b-1}}{(2b-1)!} dx = -\int_{0}^{\pi/2} \cos^{2n-1}(x) \frac{2\cdot(2x)^{2b-1}}{(2b-1)!} d(\cos x) \\
&= -\int_{0}^{\pi/2} (2n-1)\cos^{2n-1}(x) \sin(x)  \frac{2\cdot(2x)^{2b-1}}{(2b-1)!} dx + \int_{0}^{\pi/2} \cos^{2n}(x)  \frac{2^2\cdot(2x)^{2b-2}}{(2b-2)!} dx \\
&= -(2n-1) \int_{0}^{\pi/2} \cos^{2n-1}(x) \sin(x)  \frac{2\cdot(2x)^{2b-1}}{(2b-1)!} dx + 4I_{n,b-1}.
\end{align*}
Therefore,
\begin{equation}\label{Equation_Inb_2}
 \int_{0}^{\pi/2} \cos^{2n-1}(x) \sin(x) \frac{2\cdot(2x)^{2b-1}}{(2b-1)!} dx = \frac{2}{n} I_{n,b-1}.
\end{equation}
Putting \eqref{Equation_Inb_1} and \eqref{Equation_Inb_2} together, we obtain
\[ I_{n,b} = \frac{2n-1}{2n} I_{n-1,b} - \frac{1}{n^2} I_{n,b-1}, \]
or
\[ \frac{4^n}{\binom{2n}{n}} I_{n,b} = \frac{4^{n-1}}{\binom{2n-2}{n-1}}I_{n-1,b} - \frac{1}{n^2}\frac{4^n}{\binom{2n}{n}} I_{n,b-1}.  \]
By induction hypotheses, the right-hand side above is
\vspace{0.2cm}
\begin{align*}
&\frac{\pi}{2} \sum_{n-1<m_1<\cdots<m_{b}} \frac{1}{m_1^2m_2^2 \cdots m_{b}^2} - \frac{1}{n^2} \frac{\pi}{2} \sum_{n<m_2<\cdots<m_{b}} \frac{1}{m_2^2 \cdots m_{b}^2} \\
= &\frac{\pi}{2} \sum_{n<m_1<\cdots<m_{b}} \frac{1}{m_1^2m_2^2 \cdots m_{b}^2},
\end{align*}
which completes the proof of the lemma \ref{Lemma_arccos_odd_moments}.$\blacksquare$
\vspace{0.3cm}
\vspace{0.3cm}

\bt \label{Theorem_H(a,b)}
We have

 $$\displaystyle H(a, b)=\displaystyle\frac{\pi^{2b}2^{2a+3}}{(2a+2)!(2b+1)!}\int_0^{\frac{\pi}{2}}x^{2a+2}\left(1-\frac{2}{\pi}x\right)^{2b+1}\cot xdx.$$

\et
\vspace{0.4cm}

\textit{Proof.} We have seen that the Taylor series of even integer powers of  $\arcsin x$ have a connection with the multiple harmonic numbers. Formula \eqref{Equation_Taylor_arcsin_even_powers} can be rewritten as
\be \label{Equation_Taylor_arcsin_even_powers_a}
  \frac{(2\arcsin x)^{2a+2}}{(2a+2)!} = \sum_{n=1}^{\infty} \frac{4^n}{n^2 \binom{2n}{n}} x^{2n} \sum_{n_1<\cdots<n_a<n} \frac{1}{n_1^2 \cdots n_a^2} .
\ee
Equivalently, we show that
\vspace{0.2cm}
\[  \frac{(2 a+2) !(2 b+1) !}{2^{2 a+2b+3}} \cdot \frac{\pi}{2} H(a, b) = \int_{0}^{\pi / 2} x^{2 a+2}\left(\frac{\pi}{2}-x\right)^{2 b+1} \cot (x) d x. \]
\vspace{0.2cm}

Changing the variable $x \mapsto \arcsin x$ and noting that $\arccos x = \pi/2 - \arcsin x$, the assertion of the theorem is equivalent to

\begin{equation}\label{Equation_H(a,b)}
 \int_{0}^{1} \frac{(2\arcsin x)^{2a+2}}{(2a+2)!x} \frac{(2\arccos x)^{2b+1}}{(2b+1)!} dx = \frac{\pi}{2} H(a, b).
\end{equation}
\vspace{0.3cm}

By \eqref{Equation_Taylor_arcsin_even_powers_a}, we know that

\[  \frac{(2\arcsin x)^{2a+2}}{(2a+2)!x} = \sum_{n=1}^{\infty} \frac{4^n}{n^2 \binom{2n}{n}} x^{2n-1} \sum_{n_1<\cdots<n_a<n} \frac{1}{n_1^2 \cdots n_a^2} .\]
\vspace{0.3cm}

Substituting it in \eqref{Equation_H(a,b)} and integrating term by term, with the help of Lemma \ref{Lemma_arccos_odd_moments}, we obtain
\vspace{0.2cm}
\begin{align*}
&\int_{0}^{1} \frac{(2\arcsin x)^{2a+2}}{(2a+2)!x} \frac{(2\arccos x)^{2b+1}}{(2b+1)!} dx \\
=& \sum_{n=1}^{\infty} \frac{4^n}{n^2 \binom{2n}{n}} \sum_{n_1<\cdots<n_a<n} \frac{1}{n_1^2 \cdots n_a^2} \int_{0}^{1} x^{2n-1} \cdot \frac{(2\arccos x)^{2b+1}}{(2b+1)!} dx \\
=& \frac{\pi}{2} \sum_{n=1}^{\infty} \frac{1}{n^3} \left( \sum_{n_1<\cdots<n_a<n} \frac{1}{n_1^2 \cdots n_a^2} \right)\left( \sum_{n<m_1<\cdots<m_b} \frac{1}{m_1^2 \cdots m_b^2} \right) \\
=& \frac{\pi}{2} \sum_{n=1}^{\infty} \frac{1}{n^3} \sum_{n_1<\cdots < n_a < n < m_1 < \cdots < m_b} \frac{1}{n_1^2\cdots n_a^2 m_1^2 \cdots m_b^2} = \frac{\pi}{2} H(a,b).
\end{align*}$\qed{}$
\vspace{0.3cm}

Now, we relate this integral to $\hat{H}(a, b)$. This is given by the following
\vspace{0.3cm}

\bt \label{Theorem_hat_H(a,b)}
{For $a,b\geq 0$ integers}

\begin{equation*}
\int_{0}^{\pi/2} x^{2a+2}\Big(1-\frac{2x}{\pi}\Big)^{2b+1}\cot(x) \hspace{3pt} dx = \frac{(2a+2)!(2b+1)!}{2^{2a+3}\pi^{2b}}\hat{H}(a,b).
\label{eq:Habint}
\end{equation*}
\et{}

\vspace{0.4cm}

\textit{Proof.} This is equivalent with

$$\int_{0}^{1} x^{2a+2}(1-x)^{2b+1} \cot\left( \frac{\pi x}{2} \right) {\rm d} x = \frac{(2a+2)!(2b+1)!}{\pi^{2a+2b+3}} \widehat{H}(a,b).$$
\vspace{0.3cm}

Let $P(x) = x^{2a+2}(1-x)^{2b+1}$. By binomial theorem, we can rewrite $P(x)$ as $P(x) = \sum_{k=2a+2}^{2a+2b+3} (-1)^{k} \binom{2b+1}{k-2a-2} x^k$. It follows that for any $a+1 \leq k \leq a+b+1$ we have
\[ P^{(2k)}(0) = (2k)!\binom{2b+1}{2k-2a-2} = \frac{(2a+2)!(2b+1)!}{(2a+2b+3-2k)!}\binom{2k}{2a+2}. \]
Similarly, by rewriting $P(x) = -\sum_{k=2b+1}^{2a+2b+3} \binom{2a+2}{k-2b-1} (x-1)^{k} $ we find that for any $b+1 \leq k \leq a+b+1$ we have
\[ P^{(2k)}(1) = -(2k)! \binom{2a+2}{2k-2b-1} = - \frac{(2a+2)!(2b+1)!}{(2a+2b+3-2k)!} \binom{2k}{2b+1}. \]
\vspace{0.3cm}

Therefore, substituting $P(x) = x^{2a+2}(1-x)^{2b+1}$ in Lemma \ref{Lemma_Orr} we obtain

\begin{align*}
&\int_{0}^{1} x^{2a+2}(1-x)^{2b+1} \cot\left( \frac{\pi x}{2} \right) {\rm d} x \\
=& (2a+2)!(2b+1)! \cdot 2 \sum_{k=0}^{a+b+1} (-1)^{k} \left[ \binom{2k}{2a+2} - \left( 1 -  \frac{1}{2^{2k}} \right)\binom{2k}{2b+1} \right] \frac{\zeta(2k+1)}{(2a+2b+3-2k)! \pi^{2k+1}} \\
=& \frac{(2a+2)!(2b+1)!}{\pi^{2a+2b+3}} \widehat{H}(a,b).
\end{align*} 
 $\qed{}$
\vspace{0.3cm}

Now, Theorem \ref{Theorem_Zagier} follows immediately from Theorems \ref{Theorem_H(a,b)} and \ref{Theorem_hat_H(a,b)}.
\vspace{0.3cm}

\section{Zagier-type formula for $T(a, b)$}\label{Section_odd_variant}

Very recently, T. Murakami \cite{Murakami} proved an equivalent version of Theorem \ref{Theorem_odd_variant} using the same ideas as Zagier in \cite{Zagier1}. Define $\displaystyle\tilde{t}(k_{1}, \ldots, k_{r})=2^{|k|}t(k_{1}, \ldots, k_{r})$ and set

$$\displaystyle K(a, b):=\tilde{t}(\underbrace{2, 2, \ldots, 2}_{\text{$a$}}, 3, \underbrace{2, 2, \ldots, 2}_{\text{$b$}}), K(n)=\tilde{t}(\underbrace{2, 2, \ldots, 2}_{\text{$n$}})$$
then for all integers $a, b\geq 0$ we have

$$\displaystyle K(a, b)=2\sum_{k=1}^{a+b+1}(-1)^{k-1}\left[\binom{2k}{2a+1}+\left(1-\frac{1}{2^{2k}}\right)\binom{2k}{2b+1}\right]K(a+b+1-k)\zeta(2k+1).$$

As stated above, the proof is in the spirit as in \cite{Zagier1} by considering the generating functions

$$\displaystyle G(x, y)=\sum_{a, b\geq 0}(-1)^{a+b}K(a, b)x^{2a+1}y^{2b+1},$$

$$\displaystyle\hat{G}(x, y)=\sum_{a, b\geq 0}(-1)^{a+b}\hat{K}(a, b)x^{2a+1}y^{2b+1},$$
where $\hat{K}(a, b)$ denotes the right-hand side of the above formula. More exactly, Murakami showed the equality of $K(a, b)$ and $\hat{K}(a, b)$ by computing the two generating functions in a closed form and by proving that both are entire functions of exponential growth that agree at sufficiently many points so that we can force their equality. 
\vspace{0.2cm}
Again, by comparing the coefficients of like powers of $\lambda$ in the formula

$$\displaystyle\sin(\lambda\arcsin(x))=\lambda x\cdot {}_{2} F_{1}\left(\frac{1+\lambda}{2}, \frac{1-\lambda}{2}; \frac{1}{2}; x^2\right),$$
we derive the odd powers of arcsin,

\begin{equation} \label{Equation_Taylor_arcsin_odd_powers}
   \displaystyle  \displaystyle\frac{\arcsin^{2r+1}(x)}{(2r+1)!}=\sum_{n=0}^{\infty}\frac{\binom{2n}{n}}{(2n+1)4^n}\cdot x^{2n+1}\cdot\sum_{0\leq n_{1}<n_{2}<\ldots<n_{r}<n}\frac{1}{\prod_{i=1}^r(2n_{i}+1)^2}.
\end{equation}
\vspace{0.3cm}

Theorem \ref{Theorem_odd_variant} will follow from a combination of Theorem \ref{Theorem_T(a,b)}and Theorem \ref{Theorem_hat_T(a,b)}. The corresponding lemma to Lemma \ref{Lemma_arccos_odd_moments} is given by

\vspace{0.3cm}

\bl \label{Lemma_arccos_even_moments}
For any nonnegative integers $n$ and $b$, we have

\[ \int_{0}^{1} x^{2n}  \frac{\arccos^{2b+1}(x)}{(2b+1)!} dx = \frac{4^n}{(2n+1)^2 \binom{2n}{n}} \sum_{n<m_1<\cdots<m_b} \frac{1}{(2m_1+1)^2 \cdots (2m_b +1)^2} .\]
When $b=0$, the sum $\sum_{n<m_1<\cdots<m_b}$ is understood as $1$.
\el
\vspace{0.2cm}

\textit{Proof of Lemma \ref{Lemma_arccos_even_moments}.} Changing the variable $x \mapsto \cos x$ and integrating by parts, we have
\begin{align*}
\int_{0}^{1} x^{2n}  \frac{\arccos^{2b+1}(x)}{(2b+1)!} dx &=  \int_{0}^{\pi/2} \cos^{2n}(x) \sin(x)   \frac{x^{2b+1}}{(2b+1)!} dx  \\
&= \frac{1}{2n+1} \int_{0}^{\pi/2} \cos^{2n+1}(x)  \frac{x^{2b}}{(2b)!} dx.\\
\end{align*}
We denote \[ J_{n,b}:= \int_{0}^{\pi/2} \cos^{2n+1}(x)  \frac{x^{2b}}{(2b)!} dx.\]
We need to show that for $n,b \geq 0$,
\vspace{0.2cm}
\begin{equation}\label{Equation_Jnb}
J_{n,b} = \frac{4^n}{(2n+1)\binom{2n}{n}} \sum_{n<m_1<\cdots<m_b} \frac{1}{(2m_1+1)^2 \cdots (2m_b +1)^2} .
\end{equation}
\vspace{0.2cm}

We perform induction on $n+b$, but here we need to do more work for the case $n=0$. Clearly \eqref{Equation_Jnb} is true for $(n,b)=(0,0)$. Suppose that $n=0$ and $b > 0$, we have
\vspace{0.2cm}
\begin{align*}
J_{0,b} &= \int_{0}^{\pi/2} \cos(x) \frac{x^{2b}}{(2b)!} dx = \frac{\pi^{2b}}{2^{2b}(2b)!} - \int_{0}^{\pi/2} \sin(x) \frac{x^{2b-1}}{(2b-1)!} dx \\
&= \frac{\pi^{2b}}{2^{2b}(2b)!} - \int_{0}^{\pi/2} \cos(x) \frac{x^{2b-2}}{(2b-2)!} dx  = \frac{\pi^{2b}}{2^{2b}(2b)!} - J_{0,b-1}.
\end{align*}
By the well known fact $T(b):=t(\underbrace{2, 2, \ldots, 2}_{\text{$b$}}) = \frac{\pi^{2b}}{2^{2b}(2b)!}$ and the induction hypothesis for $J_{0,b-1}$, we have
\vspace{0.2cm}
\begin{align*}
J_{0,b} &= \sum_{0 \leq m_1 < \cdots < m_b} \frac{1}{(2m_1+1)^2 \cdots (2m_b +1)^2} - \sum_{0< m_2 < \cdots < m_b} \frac{1}{(2m_2+1)^2 \cdots (2m_b +1)^2} \\
&= \sum_{0 < m_1 < \cdots < m_b} \frac{1}{(2m_1+1)^2 \cdots (2m_b +1)^2}.
\end{align*}
\vspace{0.2cm}
So \eqref{Equation_Jnb} holds for $(0,b)$. For $b=0$ and any nonnegative integer $n$, \eqref{Equation_Jnb} is just Wallis' formula. Now we can assume $n,b>0$ and \eqref{Equation_Jnb} holds for smaller $n+b$. We have
\vspace{0.2cm}
\begin{align}
J_{n,b} &= \int_{0}^{\pi/2} \cos^{2n}(x) \frac{x^{2b}}{(2b)!} d(\sin x) \notag \\
&= \int_{0}^{\pi/2} 2n \cos^{2n-1}(x) \sin^2(x) dx - \int_{0}^{\pi/2} \cos^{2n}(x) \sin(x) \frac{x^{2b-1}}{(2b-1)!} dx \notag \\
&= 2n J_{n-1,b} - 2n J_{n,b} - \int_{0}^{\pi/2} \cos^{2n}(x) \sin(x) \frac{x^{2b-1}}{(2b-1)!} dx. \label{Equation_Jnb_1}
\end{align}
\vspace{0.2cm}
Also,
\begin{align}
&\int_{0}^{\pi/2} \cos^{2n}(x) \sin(x) \frac{x^{2b-1}}{(2b-1)!} dx = -\int_{0}^{\pi/2} \cos^{2n}(x)  \frac{x^{2b-1}}{(2b-1)!} d(\cos x) \notag \\
=& - \int_{0}^{\pi/2} 2n \cos^{2n}(x) \sin(x) \frac{x^{2b-1}}{(2b-1)!} dx + \int_{0}^{\pi/2} \cos^{2n+1}(x) \frac{x^{2b-2}}{(2b-2)!} dx \notag \\
=& -2n \int_{0}^{\pi/2} \cos^{2n}(x) \sin(x) \frac{x^{2b-1}}{(2b-1)!} dx + J_{n,b-1}. \label{Equation_Jnb_2}
\end{align}
Putting \eqref{Equation_Jnb_1} and \eqref{Equation_Jnb_2} together we obtain
\[ J_{n,b} = \frac{2n}{2n+1} J_{n-1,b} - \frac{1}{(2n+1)^2} J_{n,b-1}, \]
or,
\[ \frac{(2n+1)\binom{2n}{n}}{4^n} J_{n,b} = \frac{(2n-1)\binom{2n-2}{n-1}}{4^{n-1}}J_{n-1,b} - \frac{1}{(2n+1)^2} \frac{(2n+1)\binom{2n}{n}}{4^n} J_{n,b-1}. \]
By induction hypotheses, the right-hand side above is
\begin{align*}
&\sum_{n-1<m_1<\cdots<m_b} \frac{1}{(2m_1+1)^2 \cdots (2m_b +1)^2} - \frac{1}{(2n+1)^2} \sum_{n<m_2<\cdots<m_b} \frac{1}{(2m_2+1)^2 \cdots (2m_b +1)^2} \\
=& \sum_{n<m_1<\cdots<m_b} \frac{1}{(2m_1+1)^2 \cdots (2m_b +1)^2},
\end{align*}
which completes the proof of the lemma.$\blacksquare$
\vspace{0.2cm}

\bt \label{Theorem_T(a,b)}
We have

$$\displaystyle T(a,b)=\frac{\pi^{2b+1}}{2^{2b+1}(2a+1)!(2b+1)!}\int_{0}^{\pi/2} x^{2a+1}\left(1-\frac{2x}{\pi}\right)^{2b+1} \cot(x) \hspace{3pt} dx.$$

\et
\vspace{0.4cm}

\textit{Proof.} Now, in this case we will use the Taylor series expansion for odd integer powers of $\arcsin x$ and we can rewrite formula \eqref{Equation_Taylor_arcsin_odd_powers} as

\[ \frac{\arcsin^{2a+1}(x)}{(2a+1)!} = \sum_{n=0}^{\infty} \frac{\binom{2n}{n}}{(2n+1)4^n} x^{2n+1} \sum_{0\leq n_1 < \cdots < n_a < n} \frac{1}{(2n_{1}+1)^2 \cdots (2n_{a}+1)^2}. \]

Equivalently, we show that
\[ (2a+1)!(2b+1)!T(a,b) = \int_{0}^{\pi/2} x^{2a+1} \left( \frac{\pi}{2} - x \right)^{2b+1} \cot(x) dx. \]
\vspace{0.2cm}
Changing the variable $x \mapsto \arcsin x$ and noting that $\arccos x = \pi/2 - \arcsin x$, the assertion of the theorem is equivalent to

\begin{equation}\label{Equation_T(a,b)}
 \int_{0}^{1} \frac{\arcsin^{2a+1}(x)}{(2a+1)!x} \frac{\arccos^{2b+1}(x)}{(2b+1)!} dx = T(a,b).
\end{equation}

Substituting
\[ \frac{\arcsin^{2a+1}(x)}{(2a+1)!x} = \sum_{n=0}^{\infty} \frac{\binom{2n}{n}}{(2n+1)4^n} x^{2n} \sum_{0 \leq n_1 < \cdots < n_a < n} \frac{1}{(2n_1+1)^2 \cdots (2n_{a}+1)^2} \]
in \eqref{Equation_T(a,b)} and integrating term by term, with the help of Lemma \ref{Lemma_arccos_even_moments}, we have

\begin{align*}
&\int_{0}^{1} \frac{\arcsin^{2a+1}(x)}{(2a+1)!x} \frac{\arccos^{2b+1}(x)}{(2b+1)!} dx \\
=& \sum_{n=0}^{\infty} \frac{\binom{2n}{n}}{(2n+1)4^n} \sum_{0 \leq n_1 < \cdots < n_a < n} \frac{1}{(2n_1+1)^2 \cdots (2n_{a}+1)^2} \int_{0}^{1} x^{2n}  \frac{\arccos^{2b+1}(x)}{(2b+1)!} dx \\
=& \sum_{n=0}^{\infty} \frac{1}{(2n+1)^3} \left( \sum_{0 \leq n_1 < \cdots < n_a < n} \frac{1}{(2n_1+1)^2 \cdots (2n_{a}+1)^2} \right)\left( \sum_{n<m_1<\cdots<m_b} \frac{1}{(2m_1+1)^2 \cdots (2m_b +1)^2} \right) \\
=& \sum_{n=0}^{\infty} \frac{1}{(2n+1)^3} \sum_{0 \leq n_1 < \cdots < n_a < n < m_1 < \cdots < m_b} \frac{1}{(2n_1+1)^2 \cdots (2n_{a}+1)^2(2m_1+1)^2 \cdots (2m_b +1)^2}\\
=& T(a,b),
\end{align*} and the proof of Theorem \ref{Theorem_T(a,b)} is complete.$\qed{}$
\vspace{0.3cm}

Now, we relate this integral to the right-hand side quantity of Theorem \ref{Theorem_odd_variant} denoted by $\hat{T}(a,b)$. This is given by the following
\vspace{0.3cm}

\bt \label{Theorem_hat_T(a,b)}
We have

\begin{equation*}
\int_{0}^{\pi/2} x^{2a+1}\left(1-\frac{2x}{\pi}\right)^{2b+1} \cot(x) \hspace{3pt} dx = \frac{2^{2b+1}(2a+1)!(2b+1)!}{\pi^{2b+1}}\hat{T}(a,b) 
\end{equation*}
\et
\vspace{0.4cm}

\textit{Proof.} This is equivalent with

\[ \int_{0}^{1} x^{2a+1}(1-x)^{2b+1} \cot\left( \frac{\pi x}{2} \right) {\rm d} x = \frac{(2a+1)!(2b+1)!2^{2a+2b+3}}{\pi^{2a+2b+3}} \widehat{T}(a,b). \]
\vspace{0.3cm}

The proof is similar to that of Theorem \ref{Theorem_hat_H(a,b)}. Consider $P(x) = x^{2a+1}(1-x)^{2b+1}$ in Lemma \ref{Lemma_Orr} and note that
\begin{align*}
P^{(2k)}(0) &= -(2k)!\binom{2b+1}{2k-2a-1} = -\frac{(2a+1)!(2b+1)!}{(2a+2b+2-2k)!}\binom{2k}{2a+1}, \text{~for~} a+1 \leq k \leq a+b+1, \\
P^{(2k)}(1) &= -(2k)!\binom{2a+1}{2k-2b-1} = -\frac{(2a+1)!(2b+1)!}{(2a+2b+2-2k)!}\binom{2k}{2b+1}, \text{~for~} b+1 \leq k \leq a+b+1.
\end{align*}
This implies that

\begin{align*}
&\int_{0}^{1} x^{2a+1}(1-x)^{2b+1} \cot\left( \frac{\pi x}{2} \right) {\rm d} x \\
=& (2a+1)!(2b+1)! \cdot 2\sum_{k=1}^{a+b+1}(-1)^{k+1}\left[ \binom{2k}{2a+1}+\left(1 - \frac{1}{2^{2k}}\right)\binom{2k}{2b+1} \right]\frac{\zeta(2k+1)}{(2a+2b+2-2k)!\pi^{2k+1}} \\
=& \frac{(2a+1)!(2b+1)!2^{2a+2b+3}}{\pi^{2a+2b+3}} \widehat{T}(a,b).
\end{align*}
$\qed$
\vspace{0.3cm}

Now, Theorem \ref{Theorem_odd_variant} follows immediately from Theorems \ref{Theorem_T(a,b)} and \ref{Theorem_hat_T(a,b)}.

\bigskip

\section{Comments and remarks}

An interesting byproduct of our proofs for Theorem \ref{Theorem_Zagier} and Theorem \ref{Theorem_odd_variant} is some interesting arithmetic property of the numbers $\zeta(2k+1)/\pi^{2k+1}$. We conjecture that all these numbers are transcendental, but as far as things are going, very little is known. At this point, we do not even know if at least one of the values $\zeta(2k+1)/\pi^{2k+1}$ is irrational. 

Now, suppose that all of $$\displaystyle\frac{\zeta(2a+3)}{\pi^{2a+3}}, \frac{\zeta{2a+5}}{\pi^{2a+5}}, \ldots, \frac{\zeta(4a+3)}{\pi^{4a+3}}$$ are rational numbers and denote by $q_{a}$ their common denominator. By applying Lemma \ref{Lemma_Orr} for $P(x)=x^{2a+2}(1-x)^{2a+1}$ and multiplying by $2^{4a+2}q_{a}$ we obtain

\begin{equation}\label{Arithmetic_Orr}
    2^{4a+2}q_{a}\int_0^1x^{2a+2}(1-x)^{2a+1}\cot\left(\frac{\pi x}{2}\right)\in\mathbb{Z}(2a+2)!q_{a}\frac{\zeta(2a+3)}{\pi^{2a+3}}+\ldots+\mathbb{Z}(4a+2)!q_{a}\frac{\zeta(4a+3)}{\pi^{4a+3}}. 
\end{equation}

Te right-hand side of \ref{Arithmetic_Orr} is an integral multiple of $(2a+2)!$ while the left-hand side is positive and less than $q_{a}$ by employing the elementary inequality $x\cot(\pi x/2)\leq\frac{2}{\pi}<1$ for $x\in (0, 1]$. This implies that $q_{a}\geq (2a+2)!$, and the growth of the common denominators $q_{a}$ (as $a\to\infty$) is essentially the same as the one of Haynes and Zudilin in \cite{Haynes-Zudilin}. This means that the special case $H(a, a)=\hat{H}(a, a)$ of Zagier's formula already implies the above arithmetic property.

\vspace{1 em}





\bigskip

\vspace*{3mm}
\begin{flushright}
\begin{minipage}{148mm}\sc\footnotesize
Department of Mathematical Sciences, Tsinghua University, Beijing, China\\
{\it E--mail address}: {\tt lilaimath@gmail.com, lail21@mails.tsinghua.edu.cn} \vspace*{3mm}
\end{minipage}
\end{flushright}

\vspace*{3mm}
\begin{flushright}
\begin{minipage}{148mm}\sc\footnotesize
Yanqi Lake Beijing Institute of Mathematical Sciences and Applications (BIMSA) \& 
Yau Mathematical Sciences Center (YMSC), Tsinghua University, Beijing, China\\
{\it E--mail address}: {\tt lupucezar@gmail.com, lupucezar@bimsa.cn} \vspace*{3mm}
\end{minipage}
\end{flushright}

\vspace*{3mm}
\begin{flushright}
\begin{minipage}{148mm}\sc\footnotesize
University of Pittsburgh, Mathematics Department, Pittsburgh, PA, USA\\
Bank of New York Mellon, Pittsburgh, PA, USA\\
{\it E--mail address}: {\tt rrokered@gmail.com, derek.orr@bnymellon.com} \vspace*{3mm}
\end{minipage}
\end{flushright}

\end{document}